%
\def\date{8 November 2010}  
\magnification=1200
\overfullrule=0pt
\input epsf.tex
\newif\ifproofmode
\def\xrefsfilename{nonplan.xrf}  
\def\myinput#1{\immediate\openin0=#1\relax
   \ifeof0\write16{Cannot input file #1.}
   \else\closein0\input#1\fi}
\newcount\referno
\newcount\thmno
\newcount\secno
\newcount\figno
\referno=0
\thmno=0
\secno=0
\def\ifundefined#1{\expandafter\ifx\csname#1\endcsname\relax}
\myinput \xrefsfilename
\immediate\openout1=\xrefsfilename
\def\bibitem#1#2\par{\ifundefined{REFLABEL#1}\relax\else
 \global\advance\referno by 1\relax
 \immediate\write1{\noexpand\expandafter\noexpand\def
 \noexpand\csname REFLABEL#1\endcsname{\the\referno}}
 \global\expandafter\edef\csname REFLABEL#1\endcsname{\the\referno}
 \item{\the\referno.}#2\ifproofmode [#1]\fi\fi}
\def\cite#1{\ifundefined{REFLABEL#1}\ignorespaces
   \global\expandafter\edef\csname REFLABEL#1\endcsname{?}\ignorespaces
   \write16{ ***Undefined reference #1*** }\fi
 \csname REFLABEL#1\endcsname}
\def\nocite#1{\ifundefined{REFLABEL#1}\ignorespaces
   \global\expandafter\edef\csname REFLABEL#1\endcsname{?}\ignorespaces
   \write16{ ***Undefined reference #1*** }\fi}
\def\newthm#1#2\par{\global\advance\thmno by 1\relax
 \immediate\write1{\noexpand\expandafter\noexpand\def
 \noexpand\csname THMLABEL#1\endcsname{(\the\secno.\the\thmno)}}
 \global\expandafter\edef\csname THMLABEL#1\endcsname{(\the\secno.\the\thmno)}
 \bigbreak\penalty-500\noindent{\bf(\the\secno.\the\thmno)\enspace}\ignorespaces
 \ifproofmode {\bf[#1]} \fi{\sl#2}
 \medbreak\penalty-200}
\def\newdef#1#2\par{\global\advance\thmno by 1\relax
 \immediate\write1{\noexpand\expandafter\noexpand\def
 \noexpand\csname THMLABEL#1\endcsname{(\the\secno.\the\thmno)}}
 \global\expandafter\edef\csname THMLABEL#1\endcsname{(\the\secno.\the\thmno)}
 \bigbreak\penalty-500\noindent{\bf(\the\secno.\the\thmno)\enspace}\ignorespaces
 \ifproofmode {\bf[#1]} \fi{#2}
 \medbreak\penalty-200}
\def\newsection#1#2\par{\global\advance\secno by 1\relax
 \immediate\write1{\noexpand\expandafter\noexpand\def
 \noexpand\csname SECLABEL#1\endcsname{\the\secno}}
 \global\expandafter\edef\csname SECLABEL#1\endcsname{\the\secno}
 \vskip0pt plus.3\vsize
 \vskip0pt plus-.3\vsize\bigskip\bigskip\vskip\parskip\penalty-250
 \message{\the\secno. #2}\thmno=0
 \centerline{\bf\the\secno. #2\ifproofmode {\rm[#1]} \fi}
 \nobreak\smallskip\noindent}
\def\refthm#1{\ifundefined{THMLABEL#1}\ignorespaces
 \global\expandafter\edef\csname THMLABEL#1\endcsname{(?)}\ignorespaces
 \write16{ ***Undefined theorem label #1*** }\fi
 \csname THMLABEL#1\endcsname}
\def\refsec#1{\ifundefined{SECLABEL#1}\ignorespaces
 \global\expandafter\edef\csname SECLABEL#1\endcsname{(?)}\ignorespaces
 \write16{ ***Undefined section label #1*** }\fi
 \csname SECLABEL#1\endcsname}
\def\newfig#1{\global\advance\figno by 1\relax
 \immediate\write1{\noexpand\expandafter\noexpand\def
 \noexpand\csname FIGLABEL#1\endcsname{\the\figno}}\ignorespaces
 \global\expandafter\edef\csname FIGLABEL#1\endcsname{\the\figno}\ignorespaces
 \the\figno\ifproofmode{\bf[#1]}\fi}
\def\reffig#1{\ifundefined{FIGLABEL#1}\ignorespaces
 \global\expandafter\edef\csname FIGLABEL#1\endcsname{??}\ignorespaces
 \write16{ ***Undefined figure label #1*** }\fi
 \csname FIGLABEL#1\endcsname}

\def\afc{almost $4$-connected}

\def\restriction{|}
\def\he{homeomorphic embedding}
\def\emb{\hookrightarrow}

\def\junk#1{}
\font\smallrm=cmr8
\def\dfn#1{{\sl#1}}
\def\cond#1#2\par{\smallbreak\noindent\rlap{\rm(#1)}\ignorespaces
\hangindent=30pt\hskip30pt{\rm#2}\smallskip}
\def\claim#1#2\par{{\medbreak\noindent\rlap{\rm(#1)}\ignorespaces
\rightskip20pt
\hangindent=20pt\hskip20pt{\ignorespaces\sl#2}\smallskip}}
\def\proof{\smallbreak\noindent{\sl Proof. }}

\def\qed{\hfill$\square$\bigskip\medskip}
\def\sqr#1#2{{\vcenter{\vbox{\hrule height.#2pt
\hbox{\vrule width.#2pt height #1pt \kern#1pt
\vrule width.#2pt}
\hrule height.#2pt}}}}
\def\square{\mathchoice\sqr56\sqr56\sqr{2.1}3\sqr{1.5}3}
\outer\def\beginsection#1\par{\vskip0pt plus.3\vsize
   \vskip0pt plus-.3\vsize\bigskip\bigskip\vskip\parskip
   \message{#1}\centerline{\bf#1}\nobreak\smallskip\noindent}
\newcount\remarkno
\def\REMARK#1{{%
   \footnote{${}^{\the\remarkno}$}{\baselineskip=11pt #1
   \vskip-\baselineskip}\global\advance\remarkno by1}}

\nopagenumbers
\footline={\hfil}
\baselineskip=12pt
\centerline{{\bf LARGE NON-PLANAR GRAPHS AND}}
\centerline{{\bf AN APPLICATION	TO CROSSING-CRITICAL GRAPHS}}
\bigskip
\centerline{Guoli Ding%
$^{1}$\vfootnote{$^1$}{\smallrm Partially supported
by NSF under Grant No.~DMS-1001230.}}
\centerline{Department of Mathematics}
\centerline{Louisiana State University}
\centerline{Baton Rouge, Louisiana 70803, USA}
\bigskip
\centerline{Bogdan Oporowski}
\centerline{Department of Mathematics}
\centerline{Louisiana State University}
\centerline{Baton Rouge, Louisiana 70803, USA}
\bigskip
\centerline{Robin Thomas%
$^{2}$\vfootnote{$^2$}{\smallrm Partially
supported
by NSF under Grants No.~DMS-9623031 and~DMS-0701077, and by NSA under
Grant No.~MDA904-98-1-0517.
}}
\centerline{School of Mathematics}
\centerline{Georgia Institute of Technology}
\centerline{Atlanta, Georgia  30332, USA}
\bigskip
\centerline{and}
\bigskip
\centerline{Dirk Vertigan%
$^{3}$\vfootnote{$^3$}{\smallrm Partially supported
by NSA under Grant No.~MDA904-97-I-0042.
}}
\centerline{Department of Mathematics}
\centerline{Louisiana State University}
\centerline{Baton Rouge, Louisiana 70803, USA}

\vskip .3in
\centerline{\bf ABSTRACT}
\bigskip
We prove that, for every positive integer $k$, there is an integer $N$
such that every $4$-connected non-planar graph with at least 
$N$ vertices has a minor isomorphic to $K_{4,k}$,
the graph obtained from a cycle of length $2k+1$ by adding an edge
joining every pair of vertices at distance exactly $k$, or
the graph obtained from a cycle of length $k$ by adding two vertices 
adjacent to
each other and to every vertex on the cycle.
We also prove a version of this for subdivisions rather than minors,
and relax the connectivity to allow $3$-cuts with one side planar
and of bounded size.
We deduce that for every integer $k$ there are only finitely many
$3$-connected $2$-crossing-critical graphs with no subdivision isomorphic to
the graph obtained from a cycle of length $2k$ by joining all pairs of
diagonally opposite vertices.

\vfill
\noindent 5 November 1998, 
Revised \date. To appear in {\sl J.~Combin.\ Theory Ser.~B}.
\vfil\eject
\baselineskip 18pt
\footline{\hss\tenrm\folio\hss}

\newsection{intro}INTRODUCTION

In this paper \dfn{graphs} are finite  and may have loops or multiple edges.
A graph is a \dfn{subdivision} of another
if the first
can be obtained from the second by replacing each edge by a non-zero length
path
with the same ends.
Our first theorem follows the pattern of the following results.
The first two are easy.

\newthm{1.1}For every positive integer $k$, there is an integer $N$ such
that
every connected graph with at least $N$ vertices has either a path
on $k$ vertices, or a vertex with at least $k$ distinct neighbors.


\newthm{1.2}For every positive integer $k$, there is an integer $N$ such
that
every 2-connected graph with at least $N$ vertices has either a cycle
of length at least $k$, or a subgraph
isomorphic to a subdivision of  $K_{2,k}$.

These two results were generalized to $3$- and $4$-connected graphs
in~[\cite{OpoOxlTho}].
To state the theorems we need to define a few families of graphs.
Let
$k\ge 3$ be an
integer.  The \dfn{$k$-spoke wheel}, denoted by $W_k$, has
vertices $v_0,v_1,\ldots, v_k$, where $v_1,v_2,\ldots, v_k$
form a cycle, and $v_0$ is adjacent to all of
$v_1,v_2,\ldots, v_k$. 
The \dfn{$2k$-spoke alternating double wheel},
denoted by $A_k$, has vertices
$v_0,v'_0, v_1,v_2,\ldots, v_{2k}$, where $v_1,v_2,\ldots, v_{2k}$
form a cycle in this order, $v_0$ is adjacent to $v_1,v_3,\ldots, v_{2k-1}$,
and
$v'_0$ is adjacent to $v_2,v_4,\ldots,v_{2k}$.  
The vertices $v_0$ and $v'_0$ will be called the \dfn{hubs} of $A_k$.
The \dfn{$k$-rung ladder},
denoted by $L_k$, has vertices $v_1,v_2,\ldots,v_k,u_1,u_2,\ldots, u_k$,
where $v_1,v_2,\ldots,v_k$ and $u_1,u_2,\ldots,u_k$ form paths
in the order listed, and $v_i$ is adjacent to $u_i$ for $i = 1,2,\ldots,k$.
The graph
$W'_k$ is obtained from $L_k$ by adding an edge between $v_1$ and $v_k$,
and
contracting the edges joining $u_1$ to $v_1$ and $u_k$ to $v_k$. The graph
$O_k$,
called the \dfn{$k$-rung circular ladder}, is obtained from $L_k$ by adding
edges
between $v_1$ and $v_k$ and between $u_1$ and $u_k$; and the \dfn{$k$-rung
M\"obius
ladder}, denoted by $M_k$, is obtained from $L_k$ by adding  edges between
$v_1$ and $u_k$ and between $u_1$ and $v_k$. The graph $K'_{4,k}$ is
obtained from $K_{4,k}$ by splitting each of the $k$ vertices of degree
four in the same way. More precisely, it has vertices
$x,y,x',y', v_1,v_2,\ldots,v_k, v'_1,v'_2,\ldots,v'_k$, where $v_i$ is
adjacent
to $v'_i, x,$ and $y$, and  $v'_i$ is adjacent to
$v_i, x',$ and $y'$ for $i=1,2,\ldots,k$. We remark that $W_k$, $W'_k$, and
$K_{3,k}$ are 3-connected.
The following is proved in~[\cite{OpoOxlTho}].

\newthm{1.3}For every integer $k\ge 3$, there is an integer $N$ such that
every 3-connected
graph with at least $N$ vertices has a subgraph isomorphic to a subdivision
of one of $W_k$, $W'_k$,  and $K_{3,k}$.

For the second result we need a couple more definitions.
A \dfn{separation} of a graph is a pair
$(A,B)$ of subsets of $V(G)$ such that $A\cup B = V(G)$, and there is
no edge between $A-B$ and $B-A$. It is \dfn{nontrivial} if
$A-B\ne\emptyset\ne B-A$. The \dfn{order} of $(A,B)$
is $|A\cap B|$.  A graph $G$ is said to be
\dfn{almost 4-connected} if it is 3-connected and, for every separation
$(A,B)$ of $G$ of order three, one of $A-B, B-A$ contains at most one vertex.
(We remark that this is called ``internally $4$-connected" in
[\cite{OpoOxlTho}], but that term usually has a different meaning.)
Clearly every 4-connected graph is almost 4-connected,
and if $k \ge 4$, then $A_k, O_k, M_k, K_{4,k},$
and $K'_{4,k}$
are almost 4-connected.
The following is the second result from~[\cite{OpoOxlTho}].

\newthm{typical}For every integer $k\ge 4$, there is an integer $N$ such that
every almost 4-connected
graph with at least $N$ vertices contains a subgraph isomorphic to a
subdivision of one of
$A_k$, $O_k$, $M_k$, $K_{4,k}$, and $K'_{4,k}$.

Our first objective is to prove a version of~\refthm{typical} for
non-planar graphs, as follows.
We define $B_k$ to be the graph obtained from $A_k$ by adding an edge
joining its hubs.

\newthm{main1} For every integer $k\ge 4$, there is an integer $N$ such that
every almost 4-connected non-planar
graph with at least $N$ vertices has a subgraph isomorphic to a
subdivision of one of
$B_k$, $M_k$, $K_{4,k}$, and $K'_{4,k}$.

A graph is a \dfn{minor} of another if the first can be obtained from
a subgraph of the second by contracting edges.
For the minor containment~\refthm{main1} has the following corollary,
which was stated for $4$-connected graphs in the abstract.

\newthm{minorcor} For every integer $k\ge 4$, there is an integer $N$ such that
every almost 4-connected non-planar
graph with at least $N$ vertices has a minor isomorphic to $K_{4,k}$, or
the graph obtained from a cycle of length $2k+1$ by adding an edge
joining every pair of vertices at distance exactly $k$, or
the graph obtained from a cycle of length $k$ by adding two vertices 
adjacent to
each other and to every vertex on the cycle.

\proof
This follows immediately from~\refthm{main1}, because 
$K_{4,k}$ is a minor of $K'_{4,k}$;
the second outcome graph is a minor of $M_{2k+1}$;
and the third outcome graph is a minor of $B_{2k}$.~\qed
\medskip

In fact, in~\refthm{largenonplanar} we prove a stronger result 
than~\refthm{main1}.
We relax the connectivity requirement on $G$ to allow separations 
of order three as long as one side
of the separation is planar and has bounded size.

We apply the stronger form of~\refthm{main1} to deduce a theorem about
$2$-crossing-critical graphs.
Traditionally, a graph $G$ is called \dfn{$2$-crossing-critical} if it cannot
be drawn in the plane with at most one crossing, but $G\backslash e$ can
be so drawn for every edge $e\in E(G)$.
(We use $\backslash$ for deletion.
In {\sl drawings} of graphs edges are permitted to cross, whereas in
{\sl embeddings} they are not.)
But then every graph obtained from a $2$-crossing-critical graph by subdividing
an edge is again $2$-crossing-critical, and (iv) below suggests another simple 
operation that
can be used to generate arbitrarily large $2$-crossing-critical graphs.
To avoid these easily understood constructions we define
a graph $G$ to be \dfn{X-minimal} if
\item{(i)}$G$ has crossing number at least two,
\item{(ii)}$G\backslash e$ has crossing number at most one for every
edge $e\in E(G)$,
\item{(iii)}$G$ has no vertices of degree two, and
\item{(iv)}$G$ does not have a vertex of degree four
incident with two pairs of parallel edges.

\noindent
If $v$ is a vertex of degree two in a graph $G$, and $G'$ is obtained
from $G$ by contracting one of the edges incident with $v$, then $G$
satisfies (i) if and only $G'$ satisfies (i), and the same holds for
condition (ii).
Similarly, if $u\in V(G)$ has degree four and is incident with two pairs
of parallel edges, and if $G''$ is obtained from $G\backslash u$ by
adding a pair of parallel edges joining the two neighbors of $u$,
then the same conclusion holds for $G$ and $G''$.
Thus the notion of X-minimality provides a reasonable concept of 
being ``minimal with crossing number at least two".
Our second result then states the following.

\newthm{main2}
For every integer $k$ there exists an integer $N$ such that every
X-minimal graph on at least $N$ vertices has a subgraph isomorphic to 
a subdivision of $M_k$.

This is of interest, because of a belief by some experts on crossing 
numbers that X-minimal graphs with an $M_{7}$ subdivision can be completely
described. There are infinitely many of them, but they all seem to fall
within a well-described infinite family.
The sequel to~[\cite{BokOpoRicSal}] promises to prove that.
Another proof of~\refthm{main2} appears in~[\cite{BokOpoRicSal}].

To prove~\refthm{main2} we need the following lemma, which may be of
independent interest.

\newthm{main3}
Let $G$ be an {\rm X}-minimal graph on at least $17$ vertices.
Then for every separation $(A,B)$ of $G$ of order at most three,
one of $G\restriction A$, $G\restriction B$
has at most $129$ vertices and can be embedded in a disk with $A\cap B$
embedded on the boundary of the disk.

\noindent
The bound of $129$ is far from best possible, and we make no attempt
to optimize it.

The paper is organized as follows.
In Section~\refsec{planar} we state two lemmas from other papers that
will be used later.
In Section~\refsec{large} we prove~\refthm{main1}, and
in Section~\refsec{plan} we prove a lemma about planar graphs
that we use in the final Section~\refsec{xing}, where we 
first prove~\refthm{main3} and then~\refthm{main2}.

The ideas of our paper were initially developed in November 1998 and
written in manuscript form~[\cite{DinOpoThoVer}].
In October 2009 the authors of~[\cite{BokOpoRicSal}] kindly informed us
of their work, and that prompted us to revise~[\cite{DinOpoThoVer}],
resulting in the present article.


\newsection{planar}PLANAR SUBGRAPHS OF NON-PLANAR GRAPHS

We formalize the concept of a subdivision as follows.
Let $G,H$ be graphs. A mapping $\eta$ with domain $V(G)\cup E(G)$
is called a \dfn{homeomorphic
embedding} of $G$ into $H$ if for every two vertices $v,v'$  and
every two
edges $e,e'$ of $G$
\item{(i)}$\eta(v)$ is a vertex of $H$, and if $v,v'$ are distinct
then $\eta(v),\eta(v')$ are distinct,
\item{(ii)}if $e$ has ends $v,v'$, then $\eta(e)$ is a path of $H$
with ends $\eta(v),\eta(v')$, and otherwise disjoint from $\eta(V(G))$, and
\item{(iii)}if $e,e'$ are distinct, then $\eta(e)$ and $\eta(e')$ are
edge-disjoint, and if they have a vertex in common, then this vertex
is an end of both.

\noindent We shall denote the fact that $\eta$ is a
homeomorphic embedding of $G$ into $H$ by writing $\eta:G\emb H$.
If $K$ is a subgraph of $G$
we denote by $\eta(K)$ the subgraph of $H$ consisting of all
vertices $\eta(v)$, where $v\in V(K)$, and all vertices and
edges that belong
to $\eta(e)$ for some $e\in E(K)$.
It is easy to see that $H$ has a subgraph isomorphic to a
subdivision of $G$ if and only
if there is a homeomorphic embedding $G\emb H$.
The reader is advised to notice that $V(\eta(K))$ and $\eta(V(K))$
mean different sets. The first is the vertex-set of the graph
$\eta(K)$, whereas the second is the image of the vertex-set of
$K$ under the mapping $\eta$.
An \dfn{$\eta$-path} in $H$ is a path in $H$ with both ends in 
$\eta(G)$ and otherwise disjoint from it.

A cycle $C$ in a graph $G$ is called \dfn{peripheral} if it is induced
and $G\backslash V(C)$ is connected.
Let $\eta:G\emb H$,
let $C$ be a peripheral cycle in $G$, and let $P_1$ and $P_2$ be
two disjoint $\eta$-paths with ends $u_1,v_1$ and $u_2,v_2$,
respectively, such that $u_1,u_2,v_1,v_2$ belong to $V(\eta(C))$ and
occur on $\eta(C)$ in the order listed. In those circumstances we say
that the pair $P_1,P_2$ is an \dfn{$\eta$-cross}. We also say that
it is an \dfn{$\eta$-cross in $C$.} We say that $u_1,v_1,u_2,v_2$
are the \dfn{feet} of the cross. We say that the cross is \dfn{free} if
\item{(F1)}for $i=1,2$ there is no $e\in E(G)$ such that $P_i$
has both ends in $V(\eta(e))$, and
\item{(F2)}
whenever $e_1,e_2\in E(G)$ are such that all the feet of the cross
belong to $V(\eta(e_1))\cup V(\eta(e_2))$, then 
$e_1$ and $e_2$ have no end in common.

\noindent
The following is shown in~[\cite{RobSeyThoExt}].

\newthm{planarlemma}Let $G$ be an \afc\ planar graph on at least seven vertices,
let $H$ be
a non-planar graph, and let $\eta: G\emb H$ be a \he. Then there
exists a \he\ $\eta':G\emb H$  such that 
$\eta(v)=\eta'(v)$ for every vertex $v\in V(G)$ of degree at least four
and one of the following
conditions holds:
\item{(i)}there exists an $\eta'$-path in $H$ such that both of its ends
belong to $V(\eta'(C))$ for no peripheral cycle $C$ in $G$,
\item{(ii)}there exists a free $\eta'$-cross, or
\item{(iii)}there exists a separation $(X,Y)$ of $H$ of order at most
three such that
$|\eta'(V(G))\cap X-Y|\le 1$ and $H\restriction X$ does not have
an embedding in a disk with $X\cap Y$ embedded on the boundary of the disk.
\smallskip

If $\eta$ is a \he\ of $G$ into $H$, an \dfn{$\eta$-bridge}
is a connected subgraph $B$ of $H$ with
$E(B)\cap E(\eta(G)) = \emptyset$, such that either
\item{(i)}
$|E(B) | =1, E(B) = \{e\}$ say, and both ends of $e$ are in
$V(\eta (G))$, or
\item{(ii)}
for some component $C$ of $H\backslash V(\eta (G)), E(B)$
consists of all edges of $H$ with at least one end in $V(C)$.

\noindent
It follows that every edge of $H$ not in $\eta (G)$ belongs to a
unique $\eta$-bridge.
We say that a vertex $v$ of $H$ is an \dfn{attachment} of
an $\eta$-bridge $B$ if $v\in V(\eta (G)) \cap V(B)$.

Let $\eta$ be a \he\ of $G$ into $H$. We say that an $\eta$-bridge $B$
is \dfn{unstable} if there exists an edge $e\in E(G)$ such that
$V(B)\cap V(\eta(G))\subseteq V(\eta(e))$, and otherwise we say that
it is \dfn{stable}. The following result is probably due to Tutte.
A proof may be found in~[\cite{MohTho}, Lemma~6.2.1] or~[\cite{RobSeyThoExt}] 
or elsewhere.


\newthm{stable}Let $G$ be a graph, let $H$ be a simple $3$-connected graph,
and let $\eta:G\emb H$ be a \he. Then there exists a \he\ $\eta':G\emb H$
such that every $\eta'$-bridge is stable
and $\eta(v)=\eta'(v)$ for every vertex $v\in V(G)$ of degree at least three.

\newsection{large}LARGE NON-PLANAR GRAPHS

We need the following minor strengthening of \refthm{typical}.

\newthm{variation}For every two integers $k,t\ge 4$ 
there is an integer $N$ such that
every
$3$-connected
graph with at least $N$ vertices either contains a subgraph isomorphic to a
subdivision of one of
$A_k$, $O_k$, $M_k$, $K_{4,k}$, and $K'_{4,k}$, or it has 
a separation $(A,B)$ of order at most three such that $|A|\ge t$
and $|B|\ge t$.

\proof For $t=5$ this is  \refthm{typical}. For $t>5$ the result
follows by making obvious modifications to the proof of \refthm{typical}
in [\cite{OpoOxlTho}].~\qed

\newthm{hublink}Let $k\ge4$ be an integer, let $H$ be a non-planar graph,
and let $\eta:A_{2k+1}\emb H$ be a \he. Then one of the following
holds.
\item{(i)}
There exist a \he\ $\eta':A_{k}\emb H$ and an $\eta'$-path $P$ in $H$
such that $\eta'$ maps the hubs
of $A_{k}$ to the same pair of vertices $\eta$ maps the hubs of $A_{2k+1}$ to,
and the ends of  $P$ are the images of the hubs of
$A_{k}$ under $\eta'$.
\item{(ii)}There exist a \he\ $\eta':A_{2k+1}\to H$ 
and a separation $(A,B)$ of $H$ of order at most
three such that $|\eta'(V(A_{2k+1}))\cap A-B|\le 1$ and $H\restriction A$
cannot be embedded in a disk with $A\cap B$ embedded in the boundary of 
the disk.

\proof By \refthm{planarlemma} we may assume (by replacing $\eta$
by a different \he\ that maps the hubs of $A_{2k+1}$ to the same pair of
vertices of $H$ as $\eta$) that $\eta$ satisfies
(i), (ii), or (iii) of \refthm{planarlemma}. If it satisfies (iii),
then the result holds, and so we may assume that $\eta$ satisfies
\refthm{planarlemma}(i) or \refthm{planarlemma}(ii).

Assume first that $\eta$ satisfies (2.1)(i), and let $P$ be the corresponding 
$\eta$-path. Let $v_0, v_0', v_1, v_2,\ldots , v_{4k+2}$ be as in the 
definition of $A_{2k+1}$. 
If $P$ has one end in $V(\eta(v_0v_i))-\{\eta(v_i)\}$ and the other in
$V(\eta(v_0'v_j))-\{\eta(v_j)\}$ for some $i$ and $j$, 
then $A_{2k+1}\backslash\{v_0v_i,v_0'v_j\}$ has a subgraph $A$ that is 
isomorphic to a subdivision of $A_{2k-1}$. 
Let $\eta'$ be the restriction of $\eta$ to $A$ and let $P'$ be the 
$\eta(v_0)\eta(v_0')$-path in the union of $P$, $\eta(v_0v_i)$, and $\eta(v_0'v_j)$.
Then $\eta'$ and $P'$ satisfy (i). 

Thus we may assume by symmetry that both ends of $P$ are in 
$V(\eta(A_{2k+1}\backslash v_0))-\{\eta(v_0')\}$. 
In fact, we may further assume by symmetry that both ends of $P$ are in 
$V(\eta(A_{2k+1}\backslash \{v_0,v_1,v_2,\ldots,v_{2k}\}))-\{\eta(v_0')\}$. 
Since $P\cup \eta(A_{2k+1})$ is non-planar, there exist 
$i,j\in\{2k+1,2k+2,\ldots,4k+2\}$ with $|i-j|=1$ such that $P$ is 
vertex-disjoint from $\eta(Q)$, where $Q$ is the path with vertex-set 
$\{v_0,v_i,v_j,v_0'\}$. 
Let $\eta'(x)=\eta(x)$ for all vertices and edges $x$ of 
$A_{2k+1}\backslash \{v_{2k+1},v_{2k+2},\ldots,v_{4k+2}\}$. 
We define $\eta'(v_1v_{2k})$ to be the path in $H$ with ends $\eta(v_1)$ and 
$\eta(v_{2k})$ consisting of $P$ and two subpaths of 
$\eta(G)\backslash\{\eta(v_0),\eta(v_0'),\eta(v_2),\eta(v_i)\}$. 
Then $\eta' : A_k\hookrightarrow H$ and $P'=\eta(Q)$ satisfy (i).


The argument is similar when $\eta$ satisfies \refthm{planarlemma}(ii).~\qed

\newthm{circladder}Let $k\ge 1$ be an integer, and let $H$ be a 
non-planar graph such that there exists a \he\ $\eta:O_{4k}\emb H$.
Then either $H$ has a subgraph isomorphic to a subdivision of $M_k$, or
there exist a \he\ $\eta':O_{4k}\emb H$ and
a separation $(A,B)$ of $H$ of order at most three such that
$|\eta'(V(O_{4k}))\cap A-B|\le 1$ and $H\restriction A$
cannot be embedded in a disk with $A\cap B$ embedded in the boundary 
of the disk.

\proof The proof is similar to that of \refthm{hublink}. We
omit the details.~\qed

Let us recall that $B_k$ is the graph obtained from $A_k$ by adding an edge
joining its hubs.
A graph $G$ is \dfn{$t$-shallow} if for every separation $(A,B)$ of
order at most three, one of $G\restriction A$, $G\restriction B$
has fewer than $t$ vertices and can be embedded in a disk with $A\cap B$
embedded on the boundary of the disk.
The following is the main result of this section.
It implies~\refthm{main1}, because every almost $4$-connected graph
is $5$-shallow.

\newthm{largenonplanar}For every two integers $k,t\ge 4$
there is an integer $N$ such that
every
$3$-connected $t$-shallow non-planar
graph with at least $N$ vertices contains a subgraph isomorphic to a
subdivision of one of
$B_k$, $M_k$, $K_{4,k}$, and $K'_{4,k}$.

\proof Let $k,t$ be given. By replacing $k$ by a larger integer we may
assume that $8k\ge t+1$.
Let $N$ be the integer that satisfies
\refthm{variation} with $k$ replaced by $4k$. We claim that $N$
satisfies the conclusion of \refthm{largenonplanar}. To see this
let $G$ be a $3$-connected $t$-shallow non-planar graph on at least
$N$ vertices. By \refthm{variation} $G$ has a subgraph isomorphic to
a subdivision of one of 
$A_{4k}$, $O_{4k}$, $M_{4k}$, $K_{4,4k}$, and $K'_{4,4k}$.
If $G$ has a subgraph isomorphic to
a subdivision of 
$M_{4k}$, $K_{4,4k}$, or $K'_{4,4k}$, then the
result holds.

Assume now that there exists a \he\ $\eta:A_{4k}\emb G$.
By \refthm{hublink} either $G$ has a subgraph isomorphic to
a subdivision of $B_{k}$, or there exists a separation $(A,B)$
as in \refthm{hublink}(ii). In the former case the theorem holds, 
and so we may assume the latter. Since $G$ is $t$-shallow we see
that $|B|< t$. However, all but possibly one vertex of $\eta(V(A_{4k}))$
belong to $B$,  contrary to $8k\ge t+1$. The argument when
there exists a \he\ $\eta:O_{4k}\emb G$ is similar, using \refthm{circladder}
instead.~\qed

\newsection{plan}A LEMMA ABOUT PLANAR GRAPHS

The objective of this section is to prove~\refthm{planlemma}.
Let $G$ be a {\sl plane} graph; that is, a graph embedded in the plane.
Then every cycle $C$ bounds a disk in the plane, and we define
ins$(C)$ to be the set of edges of $G$ embedded in the open disk bounded by $C$.
(By definition, an edge of an embedding or drawing does not include its ends.)
The following will be a hypothesis common to several lemmas,
and so we give it a name in order to avoid repetition.

\newdef{hypot} {\bf Hypothesis.}
Let $G$ be a loopless plane graph embedded in the closed unit disk $\Delta$,
let $x_1,x_2,x_3$ be distinct vertices of $G$, and let them be 
the only  vertices of $G$ embedded in the boundary of $\Delta$.
Assume that there is no separation $(A,B)$ of order at most two
with $x_1,x_2,x_3\in A$ and $B-A\ne\emptyset$.

\noindent
The last assumption of~\refthm{hypot} will be referred to as
the \dfn{internal $3$-connectivity of $G$}.

Assume~\refthm{hypot}, let $C$ be a cycle in $G$ with
$\{x_1,x_2,x_3\}\not\subseteq V(C)$ and ins$(C)\ne\emptyset$.
We say that $C$ is \dfn{robust} if there exists an edge $f\in\hbox{ins}(C)$
such that for every $e\in E(C)$ the graph 
$G\backslash\{x_1,x_2,x_3\}\backslash e\backslash f$
has a component containing a neighbor of each  of $x_1,x_2,x_3$.
Let $Z$ be the set of all vertices $v\in V(C)$ such that either
$v\in\{x_1,x_2,x_3\}$ or $v$ is incident with an edge not
in $E(C)\cup\hbox{ins}(C)$.
We say that $C$ is \dfn{flexible} if $|Z|\le3$ and at least two vertices in 
$Z-\{x_1,x_2,,x_3\}$ are incident with exactly one edge not in
 $E(C)\cup\hbox{ins}(C)$.
Our objective in this section is to prove that if $G$ has
sufficiently many vertices and satisfies Hypothesis~\refthm{hypot},
then it has a robust cycle or a flexible cycle.

\newthm{cyclefacial} Assume~{\rm\refthm{hypot}}.
Then every cycle of $G\backslash \{x_1,x_2,x_3\}$ that does not bound a face 
is robust.

\proof Let $C$ be a cycle of $G\backslash \{x_1,x_2,x_3\}$ that does
not bound a face, and let $f\in \hbox{ins}(C)$.
By the internal $3$-connectivity of $G$ there exist three internally disjoint paths
from $\{x_1,x_2,x_3\}$ to $V(C)$, and hence 
$G\backslash \{x_1,x_2,x_3\}\backslash e\backslash f$
has a component containing neighbors of all of $x_1,x_2,x_3$ for all $e\in E(C)$.
Thus $C$ is robust, as desired.~\qed

Let us recall that a \dfn{block} is a graph with no cut-vertices, and a
\dfn{block of a graph} is a maximal subgraph that is a block.
The \dfn{block graph} of a graph $G$ is the graph whose vertices
are all the blocks of $G$ and all the cut vertices of $G$,
with the obvious incidences.
An \dfn{end-block} of a graph $G$
is a block that has degree one in the block graph of $G$.

\newthm{cyclesdisjoint}  Assume~{\rm\refthm{hypot}}, and that
$G$ has no robust cycle.
Then every two distinct cycles of $G\backslash \{x_1,x_2,x_3\}$ 
are edge-disjoint.
Consequently, every block of $G\backslash  \{x_1,x_2,x_3\}$ is a cycle or
a complete graph on at most two vertices.

\proof This follows from~\refthm{cyclefacial}, because otherwise
some cycle of $G\backslash \{x_1,x_2,x_3\}$ is not facial.~\qed



\newthm{end-blocks} Assume~{\rm\refthm{hypot}}, and assume
that $G$ has at least $16$ vertices and no robust cycle.
Let $B_1,B_2$ be two distinct end-blocks of 
$G\backslash \{x_1,x_2,x_3\}$. For
$i=1,2$ let $v_i$ be the unique cut vertex of $G\backslash  \{x_1,x_2,x_3\}$
that belongs to $B_i$, and let
$N_i\subseteq  \{x_1,x_2,x_3\}$ be the set of vertices of $ \{x_1,x_2,x_3\}$
that have a neighbor in $B_i\backslash v_i$. 
Then $|N_1|=|N_2|=2$ and $|N_1\cap N_2|=1$.

\proof
We first notice that $N_1$ and $N_2$ have at least two
elements by the internal $3$-connectivity of $G$. 
Thus it suffices to show that $|N_1\cap N_2|\le1$. 
Let us assume for a contradiction that $x_1,x_2\in N_1\cap N_2$.
The fact that $G$ is embedded in a disk with $x_1,x_2,x_3$ on the boundary
of the disk implies that either $x_3$ has no neighbor outside 
$B_1\backslash v_1$,
or it has no neighbor outside $B_2\backslash v_2$, and hence from the symmetry
we may assume the latter.
But $x_3$ has at least one neighbor in $B_2\backslash v_2$ by the
internal $3$-connectivity of $G$.
Since $G$ has at least $16$ vertices, it follows from~\refthm{cyclesdisjoint}
that $G\backslash \{x_1,x_2,x_3\}$ has at least seven
vertices with at most two neighbors.
Each of those vertices has a neighbor in $\{x_1,x_2,x_3\}$, and hence
there is an index $i\in\{1,2,3\}$ such that
$x_i$ has at least three neighbors in 
$G\backslash \{x_1,x_2,x_3\}$.
Furthermore, 
if $B_2$ has a unique edge, then $i$ and the three neighbors of $x_i$
can be chosen to be not in $B_2\backslash v_2$.
Thus there is a cycle $C$ of $G$
containing $x_i$ but no other $x_j$
such that ins$(C)$ includes an edge $f$ incident with $x_i$;
and if $B_2$ has a unique edge, then $C$ does not use that edge.
Since $x_1,x_2$ and $x_3$ all have a neighbor in $B_2\backslash v_2$,
it follows that $C$ is robust, a contradiction.~\qed

\newthm{blockpath} Assume~{\rm\refthm{hypot}}, and assume
that $G$ has at least $16$ vertices and no robust cycle.
Then the block graph of $G\backslash \{x_1,x_2,x_3\}$ is a path.

\proof
Suppose for a contradiction that the block graph
of $G\backslash \{x_1,x_2,x_3\}$ is not a path. 
Then $G\backslash \{x_1,x_2,x_3\}$ has at least
three end-blocks, say $B_1$, $B_2$, and $B_3$. For $i=1,2,3$ let $N_i$
be as in~\refthm{end-blocks}. By~\refthm{end-blocks} we may assume
that the blocks $B_1,B_2,B_3$ are numbered in such a way that
$N_1=\{x_2,x_3\}$, $N_2=\{x_1,x_3\}$, and $N_3=\{x_1,x_2\}$.
Let $C$ be a cycle containing an edge joining $x_i$ to a vertex of
$N_j$ for all distinct integers $i,j\in\{1,2,3\}$, such that all
other edges of $C$ belong to $B_1\cup B_2\cup B_3$. 
Let $T$ be a connected subgraph of $G\backslash \{x_1,x_2,x_3\}$
such that $V(T\cap C)=\{u_1,u_2,u_3\}$, where $u_i\in V(B_i)$.
Then $x_1,u_3,x_2,u_1,x_3,u_2$ appear on $C$ in the order listed.
Since $G$ has at least $16$ vertices there exist an edge 
$f\in E(G)-E(T)-E(C)$ and index $i\in\{1,2,3\}$ such that 
$f\in \hbox{ins}(C')$, where $C'$ is the
unique cycle in $(C\cup T)\backslash u_i$.
It follows that $C'$ is robust, a contradiction.~\qed

\newthm{planlemma} Assume~{\rm\refthm{hypot}}, and assume
that $G$ has at least $130$ vertices.
Then $G$ has a  robust cycle or a flexible cycle.

\proof Assume for a contradiction that $G$ has neither a robust cycle
nor a flexible cycle.
Let $B:=G\backslash \{x_1,x_2,x_3\}$,
let $a_1b_1,a_2b_2,\ldots,a_tb_t$ be all the cut edges of $B$,
and let $D_0,D_1,\ldots,D_t$ be all the components of 
$B\backslash \{a_1b_1,a_2b_2,\ldots,a_tb_t\}$.
By~\refthm{blockpath} the numbering can be chosen so that $a_j\in V(D_{j-1})$
and $b_j\in V(D_j)$ for all $j=1,2,\ldots,t$.
By~\refthm{end-blocks} we may assume that $x_1$ and $x_3$ have a neighbor
in $D_0$, and that $x_2$ and $x_3$ have a neighbor in $D_t$.

\claim{1}For $i\in\{1,2,3\}$ and $j\in \{0,1,\ldots,t\}$ there are
at most two edges with one end $x_i$ and the other end in $D_j$.

To prove (1) suppose for a contradiction that there are three edges
with one end $x_i$ and the other end in $D_j$.
Then there exists a cycle $C$ using two of those edges such that
the third edge, say $f$, belongs to ins$(C)$ and $C\backslash x_i$ is
a subgraph of $D_j$.
If $0<j<t$, then there exists a path $P$ in $D_j\backslash E(C)$ with
ends $b_j$ and $a_{j+1}$. 
By considering the edge $f$ and path $P$ (when $0<j<t$) we deduce
that $C$ is robust, a contradiction.
This proves (1).

\claim{2}For $j=0,1,\ldots,t$ the graph $D_j$ has at most $18$ vertices.

To prove (2) we first notice that the block graph of $D_j$ is a path
by~\refthm{blockpath}. 
Since $D_j$ is $2$-edge-connected, each block of $D_j$ is a cycle
by~\refthm{cyclesdisjoint}.
By the internal $3$-connectivity of $G$
no two consecutive blocks of $D_j$ are both a cycle of length two,
unless their shared vertex is adjacent to at least one of $x_1,x_2,x_3$.
Since every vertex of $D_j$ except possibly $b_j$ (if $j>0$) and 
$a_{j+1}$ (if $j<t$)
has at least three distinct neighbors by the internal $3$-connectivity of $G$, 
the claim follows from (1).
This proves (2).

\claim{3}There is at most one index $j\in\{1,2,\ldots,t-1\}$ such that
the graph $D_j$ includes a neighbor of $x_1$.

To prove (3) we suppose for a contradiction that there exist two such
indices $j,j'$ with $0<j'<j<t$. 
Since $x_1$ has also a neighbor in $B_0$,
there exists a cycle $C$ through the vertex $x_1$ with
$V(C)\subseteq V(D_0\cup D_1\cup\cdots\cup D_j)\cup\{x_1\}$ and such that
some edge $f$ incident with $x_1$ belongs to ins$(C)$.
Since $x_2$ and $x_3$ have a neighbor in $D_t$, and $D_j$ is $2$-edge-connected,
it follows that $C$ is robust, a contradiction.
This proves (3).

 From the symmetry between $x_1$ and $x_2$ we deduce

\claim{4}There is at most one index $j\in\{1,2,\ldots,t-1\}$ such that
the graph $D_j$ includes a neighbor of $x_2$.

We are now ready to complete the proof of the lemma.
Since $G$ has at least $130$ vertices, it follows from (2)
that $t\ge8$, and hence by (3) and (4)
there exists an integer $j\in\{1,2,\ldots,t-2\}$ such that both
$D_j$ and $D_{j+1}$ include no neighbor of $x_1$ or $x_2$.
Thus each of them includes a neighbor of $x_3$ by the internal
$3$-connectivity of $G$, and hence there exists
a cycle $C$ with $V(C)\subseteq V(D_j\cup D_{j+1})\cup\{x_3\}$,
$x_3,b_j,a_{j+2}\in V(C)$, and such that $a_jb_j$ is the only edge
of $G$ incident with $b_j$ that does not belong to $E(C)\cup\hbox{ins}(C)$,
and $a_{j+2}b_{j+2}$ is the only such edge incident with $a_{j+2}$.
By considering the set $Z=\{a_{j+2},b_j,x_3\}$ we deduce
that $C$ is flexible, as desired.~\qed

We also need the following mild strengthening of~\refthm{planlemma}.
If $C$ is a subgraph of a graph $G$, then by a \dfn{$C$-bridge}
we mean an $\eta$-bridge, where $\eta:C\emb G$ is the \he\ that maps
every vertex and edge of $C$ onto itself.

\newthm{stablebridges}
Assume~{\rm\refthm{hypot}}, and let $C$ be a robust or flexible cycle
in $G$ with ${\rm ins}(C)$ maximal.
Then for every $C$-bridge $B$ of $G$ either $E(B)\subseteq{\rm ins}(C)$,
or at least one of $x_1,x_2,x_3$ belongs to $V(B)-V(C)$.

\proof
Assume first that $C$ is robust,  let $f\in \hbox{ins}(C)$ be
as in the definition of robust, and suppose for a contradiction
that $B$ is a $C$-bridge that satisfies neither conclusion of the lemma.
By the internal $3$-connectivity of $G$ the bridge $B$ includes a path
$P$ of $G\backslash\{x_1,x_2,x_3\}$
with both ends on $C$, and otherwise disjoint from it.
The graph $C\cup P$ includes a cycle $C'\ne C$ with 
ins$(C)$ properly contained in ins$(C')$.
Since every edge of $P$ belongs to a cycle of $G\backslash f$
it follows that $C'$ is robust, contrary to the maximality of $C$.

The argument when $C$ is flexible is similar. In that case the set $Z$
from the definition of flexible is the same for $C$ and $C'$.~\qed

\newsection{xing}LARGE GRAPHS WITH CROSSING NUMBER AT LEAST TWO

Recall that a graph $G$ is \dfn{X-minimal} if 
\item{(i)}$G$ has crossing number at least two,
\item{(ii)}$G\backslash e$ has crossing number at most one for every
edge $e\in E(G)$,
\item{(iii)}$G$ has no vertices of degree two, and
\item{(iv)}$G$ does not have a vertex of degree four
incident with two pairs of parallel edges.

\newthm{3conn}Every X-minimal graph on at least $17$
vertices is $3$-connected.

\proof Let $G$ be an X-minimal graph on at least $17$ vertices, 
and suppose for a contradiction that it is not $3$-connected. Thus
it has a nontrivial separation $(A,B)$ of order at most two.
We may assume that $(A,B)$ has the minimum order among all nontrivial
separations of $G$.

Assume first that the order of $(A,B)$ is at most one. Both
$G|A$ and $G|B$ have crossing number at most one by the X-minimality of
$G$. They are both non-planar, for otherwise $G$ itself would have crossing
number at most one. Thus both $G|A$ and $G|B$ have subgraphs isomorphic
to subdivisions of $K_5$ or $K_{3,3}$ by Kuratowski's theorem. 
Now the X-minimality of $G$
implies that $G|A$ and $G|B$ have at most seven vertices,
contrary to the fact that $G$ has at least $17$ vertices.

We may therefore assume that $G$ is $2$-connected and that
the order of $(A,B)$ is two. Let $A\cap B=\{u,v\}$.
Let $G_1$ be the graph obtained from $G|A$ as follows.
If $G|B$ has two edge-disjoint paths with ends $u$ and $v$, then
$G_1$ is obtained from $G|A$ by adding two edges with ends $u$ and $v$;
otherwise $G_1$ is obtained from $G|A$ by adding one edge with ends
$u$ and $v$. We define $G_2$ analogously (with the roles of $A$ and $B$
interchanged).

\claim{1}The graphs $G_1$ and $G_2$ have crossing number at most one.

To prove (1) it suffices to argue for $G_1$. 
Assume first that $G|B$ does not have two edge-disjoint paths with ends $u$ and $v$.
Since $G|B$ has a path with ends $u$ and $v$ by the $2$-connectivity of $G$, we 
deduce that a subdivision of $G_1$ is isomorphic to a subgraph of $G$,
and that the containment is proper. Thus $G_1$ has crossing number at most one
by the X-minimality of $G$.
We may therefore assume that $G|B$ has two edge-disjoint paths $P_1$ and $P_2$
with ends $u$ and $v$.
Then by choosing the paths with $P_1\cup P_2$ minimum it can be arranged that both
$P_1$ and $P_2$ pass through the vertices of $V(P_1)\cap V(P_2)$
in the same order. 
The graph $(G|A)\cup P_1\cup P_2$ is a proper subgraph of $G$ by the
X-minimality of $G$, and hence has crossing number at most one.
It follows that $G_1$ has crossing number at most one.
This proves (1).

\claim{2}The graphs $G_1$ and $G_2$ are non-planar.

To prove (2) it again suffices to argue for $G_1$. Suppose for 
a contradiction that $G_1$ is planar. By (1) there exists a planar drawing
of $G_2$ with at most one crossing. If none of the edges of $E(G_2)-E(G|B)$
is involved in the crossing, then this drawing and a planar
embedding of $G_1$ can be combined to produce a planar drawing of $G$
with at most one crossing. Thus we may assume that an edge of
$E(G_2)-E(G|B)$ is crossed by another. Therefore we may assume
that $E(G_2)-E(G|B)$ consists
of a unique edge, say $e$, and hence, by construction, $G_1$ does not have
two edge-disjoint paths with ends $u$ and $v$. By Menger's theorem
$G_1$ has an edge $f$ such that $G_1\backslash f$ has no path between
$u$ and $v$. Using the drawings of $G_1$ and $G_2$ it is now possible
to obtain a drawing of $G$, where $e$ and $f$ are the only two edges
that cross, contrary to the fact that $G$ has crossing number at least
two. This proves (2).

>From (2) and Kuratowski's theorem it follows that for $i=1,2$ the
graph $G_i$ has 
a subgraph $H_i$ isomorphic to a subdivision of $K_5$ or $K_{3,3}$.
But $H_1\cup H_2$ has crossing number at least two, and hence 
the X-minimality of $G$ implies that both $G_1$ and $G_2$ have
at most eight vertices, 
contrary to the fact that $G$ has at least $17$ vertices.
This proves that $G$ is $3$-connected.~\qed

\newthm{orangebox}Let $G$ be a graph, let $C$ be a cycle in $G$,
and let $B_0,B_1,\ldots,B_k$ be the $C$-bridges of $G$ such that
the graph $C\cup B_1\cup B_2\cup\ldots\cup B_k$ has a planar drawing with
no crossings in which
$C$ bounds a face. Let $H$ denote the graph $C\cup B_0$, and let
$f\in E(B_1)$. Assume
further
that either $G\backslash e\backslash f$ is non-planar for every $e\in E(C)$,
or that the $C$-bridge $B_0$ has exactly three attachments, two
of which have degree three in $H$. 
If $G\backslash f$ has crossing number at most one, 
then so does $G$.

\proof 
Let $\Gamma$ be a drawing of $G\backslash f$ with at most one crossing.
Our first objective is to modify $\Gamma$ to produce a drawing
of $H$ with at most one crossing such that no edge of $C$ is crossed by
another edge.
If no edge of $C$ is crossed by another edge in the drawing $\Gamma$,
then its restriction to $H$ is as desired.
Thus we may assume that an edge $e\in E(C)$ is crossed by another
edge $e'$ in $\Gamma$. It follows that $G\backslash e\backslash f$ is
planar, and hence, by hypothesis, the $C$-bridge $B_0$ has exactly 
three attachments, say $v_1,v_2,v_3$, such that $v_1$ and $v_2$
have degree three in $H$.
If $e'\not\in E(B_0)$, then it is easy to convert $\Gamma$ to
a desired drawing of $H$.
Thus we may assume that $e'\in E(B_0)$.
It follows that $B_0\backslash e'$ has two components, say $J_1$ and $J_2$,
such that $J_1$ is drawn in the closed disk bounded by $C$ and
$J_2$ is drawn in the closure of the other face of $C$.
Using the fact that $v_1$ and $v_2$ have degree three in $H$ it is now
easy to draw $J_2$ in the closed disk bounded by $C$ so as to obtain
a desired drawing of $H$.
This proves our claim that $H$ has a drawing with at most one
crossing such that no edge of $C$ is crossed by another edge in
that drawing. Thus $C$ bounds a face.
By hypothesis it is possible to draw
$B_1\cup B_2\cup\ldots\cup B_k$ without crossings in that face, 
showing that $G$ has crossing number
at most one, as desired.~\qed

Let $G$ be a graph, let $u,u_1,u_2,u_3$ be distinct vertices of $G$, 
and let $Q_1,Q_2,Q_3$ be three
paths in $G$ such that $Q_i$ has ends $u$ and $u_i$ and such that
$Q_1,Q_2,Q_3$ are disjoint except for $u$. We say that
$Q_1\cup Q_2\cup Q_3$ is a \dfn{triad} in $G$, and that
the vertices $u_1,u_2,u_3$ are its \dfn{feet}.
Let $G$ be a graph, and let $P_1,P_2,P_3$ be three pairwise disjoint 
paths in $G$,
where $P_i$ has ends $u_i$ and $v_i$.
Let $T_1$ and $T_2$ be two triads with feet $v_1,v_2,v_3$ such that
the graphs $P_1\cup P_2\cup P_3$, $T_1$, $T_2$ are pairwise disjoint,
except for $v_1,v_2,v_3$. In those circumstances we say that
$P_1\cup P_2\cup P_3\cup T_1\cup T_2$ is a \dfn{tripod}, and that
the vertices $u_1,u_2,u_3$ are its \dfn{feet}.
We need the following result of [\cite{RobSeyGM9}].

\newthm{tripod}Let $G$ be a graph, and let $u_1,u_2,u_3$
be three vertices of $G$ such that there is no separation $(A,B)$
of $G$ of order at most two with $u_1,u_2,u_3\in A$ and $B-A\ne\emptyset$.
If $G$ has no planar embedding with the
vertices $u_1,u_2,u_3$ incident with the same face, then
$G$ has a tripod with feet $u_1,u_2,u_3$.

\newthm{flat}Let $G$ be an X-minimal graph on at least $17$ vertices,
and let $(A,B)$ be a separation in $G$ of order three. Then one of
$G|A$, $G|B$ has a planar embedding with the vertices $A\cap B$ embedded
on the boundary of the same face.

\proof Suppose for a contradiction that the conclusion does not hold.
By \refthm{3conn} the graph $G$ is $3$-connected.
By \refthm{tripod} $G|A$ has a tripod $T_1$ with feet $A\cap B$,
and $G|B$ has a tripod $T_2$ with feet $A\cap B$. The graph $T_1\cup T_2$
has crossing number at least two, as is easily seen. Thus $G=T_1\cup T_2$
by the X-minimality of $G$. Moreover, the X-minimality of $G$
implies that $G$ has at most $10$ vertices, a contradiction.~\qed

We are now ready to prove~\refthm{main3}, which we restate.

\newthm{shallow}Every X-minimal graph on at least $17$
vertices is $130$-shallow.

\proof 
Let $G$ be an X-minimal graph on at least $17$ vertices,
and let $(A,B)$ be a separation
in $G$ of order at most three with $A-B\ne\emptyset\ne B-A$.
By \refthm{3conn} the separation $(A,B)$ has order exactly three.
By \refthm{flat} we may assume that $G|B$ is embedded in a disk with
the vertices of $A\cap B$ embedded in the boundary of the disk.
It follows that $G|B$ satisfies~\refthm{hypot}, where $A\cap B=\{x_1,x_2,x_3\}$.
We may and shall assume for a contradiction that $|B|\ge130$.
By~\refthm{planlemma} applied to the graph $G|B$ we deduce that
$G|B$ has a cycle $C$ that is robust or flexible.
By~\refthm{stablebridges} we may choose $C$ so
that exactly one $C$-bridge $B_0$ of $G$ satisfies 
$E(B_0)\not\subseteq\hbox{ins}(C)$.
We wish to apply~\refthm{orangebox}, and so we need to verify the hypotheses.
If $C$ is robust, then let $f$ be as in the definition of robust;
otherwise let $f\in\hbox{ins}(C)$ be arbitrary.
If $C$ is flexible, then the bridge $B_0$ has exactly three attachments,
and two of them have degree three in $C\cup B_0$.
Now let $C$ be robust, and let $e\in E(C)$.
We claim that $G\backslash e\backslash f$ is not planar.
To prove this we first notice that $G|A$ cannot be embedded in a disk with
$A\cap B$ embedded in the boundary of the disk, because $G|B$ can be so 
embedded and $G$ is not planar.
By~\refthm{tripod} the graph $G|A$ has a tripod $T$ with feet $A\cap B$.
Since $C$ is robust the graph $(G|B)\backslash e\backslash f$ has
a connected subgraph $R$ that includes $A\cap B$.
It follows that $T\cup R$ is a subdivision of $K_{3,3}$, which proves our
claim that $G\backslash e\backslash f$ is not planar.
The graph $G\backslash f$ has crossing number at most one by the
X-minimality of $G$, and hence
by~\refthm{orangebox} the graph $G$ has crossing number at most one,
a contradiction.~\qed

\newthm{Asubdiv}Let $G$ be the graph obtained from $A_4$ by subdividing
the edges $v_1v_2$ and $v_5v_6$, and joining the new vertices by
an edge. Then $G$ has crossing number at least two.

\proof This follows from the fact that the new edge is the only edge
$e\in E(G)$ such that $G\backslash e$ is planar.~\qed

\newthm{noB}No X-minimal graph has a subgraph isomorphic to a subdivision
of $B_{65}$.


\proof
Let $H$ be an X-minimal graph, and suppose for a contradiction
that it has a subgraph isomorphic to a subdivision
of $B_{65}$. Let $\eta:B_{65}\emb H$ be a \he, and let $\eta_0$
be the restriction of $\eta$ to $A_{65}$.
Let $e_0$ be the edge of $B_{65}$ joining the two hubs.
Let $J$ be the union of $\eta_0(A_{65})$ and all $\eta_0$-bridges
except the one that includes $\eta(e_0)$. We claim that
$J$ is planar. To prove this claim suppose for a contradiction that
it is not. 
 From \refthm{hublink}
applied to $A_{65}, J$, and $\eta_0$ we deduce that (i) or (ii)
of \refthm{hublink} holds. If (i) holds, then we conclude that
the graph obtained from $B_{32}$ by adding an edge parallel to $e_0$
is isomorphic to a subdivision of $H$. That is a contradiction,
because said graph is not X-minimal, as is easily seen. Thus we may
assume that~\refthm{hublink}(ii) holds; that is,
$H$ has a separation $(A,B)$ as in \refthm{hublink}(ii).
But $|B|\ge |V(B_{65})|-1\ge 130$, and $H|A$ does not have a planar 
embedding with
the vertices in $A\cap B$ incident with the same face, contrary
to \refthm{shallow}. This proves our claim that $J$ is planar.
Thus we may regard $J$ as a graph embedded in the sphere.

Let the vertices of $A_{65}$ be numbered as in the definition of $A_{65}$.
Assume first that $\eta(e_0)$ has only one edge. 
Let $C_0$ be a cycle in $J$ with
$v_0\not\in V(C_0)$ such that the open disk bounded by $C_0$
that includes $v_0$ is as small as possible. Let $C_0'$ be
defined analogously, with $v_0'$ replacing $v_0$. The cycles
$C_0,C_0'$ are edge-disjoint, for otherwise $H$ has crossing number
at most one. But now it follows that the graph obtained from $H$
by deleting an edge of $\eta(v_0v_1)$ has crossing number at least
two, contrary to the X-minimality of $H$. This completes the
case when $\eta(e_0)$ has only one edge.

We may therefore assume that $\eta(e_0)$ has at least one internal vertex.
Let us say that an $\eta$-bridge of $H$ is \dfn{solid} if either it has
at least two edges, or it has a unique edge and that edge is not parallel to an
edge of $\eta(B_{65})$.
By \refthm{stable} we may assume that every solid $\eta$-bridge is stable.
Let us say that a vertex $v\in V(\eta_0(A_{65}))-\{\eta_0(v_0),\eta_0(v_0')\}$ 
is \dfn{exposed} if there exists an $\eta$-path between an internal vertex of
$\eta(e_0)$ and $v$.
It follows from \refthm{3conn} that there exists at  least one exposed vertex.
For an integer $i\in\{1,3,\ldots,129\}$ let $C_i$ denote the cycle
of $A_{65}$ with vertex-set $\{v_i,v_{i+1},v_{i+2},v_{i+3},v_{i+4},v_0\}$
(index arithmetic modulo $130$), 
and let $F_i$ be the set of edges of $A_{65}$ with at least one end in $V(C_i)$.
>From \refthm{Asubdiv} we deduce that there exists an integer $i$
such that $\eta(e)$ includes an exposed vertex for no $e\in F_i$.
Let $J_0,J_1,\ldots,J_k$ be all the $\eta(C_i)$-bridges of $H$, where
$J_0$ includes $v_0'$.
Then $J_0$ includes $\eta(e_0)$, and hence $J_1,J_2,\ldots,J_k$ are
also $\eta_0(C_i)$-bridges of $J$.
Since every solid $\eta$-bridge is stable, it follows that $J_1,J_2,\ldots,J_k$,
when regarded as $\eta_0(C_i)$-bridges of $J$, are embedded in 
the closed disk $\Delta$ bounded by $\eta_0(C_i)$ that does not include $v_0'$;
hence $\eta_0(C_i)\cup J_1\cup J_2\cup\cdots\cup J_k$ has a planar 
embedding in which $\eta_0(C_i)$ bounds a face.
Since in the planar embedding of $J$ the path $\eta(v_0v_{i+2})$ is embedded
in $\Delta$ we deduce that $k\ge1$.
Thus we may select
$f\in E(J_1)$.
Since there exists an exposed vertex, but none in $\eta(e)$ for any
$e\in F_i$, it follows that $H\backslash e\backslash f$ is non-planar for every
edge $e\in E(C_i)$.
The graph $H\backslash f$ has
crossing number at most one by the X-minimality of $G$, contrary to
\refthm{orangebox}.~\qed

We are finally ready to prove~\refthm{main2}, which we restate.

\newthm{Msubd}For every integer $k$ there exists an integer $N$ such that
every  X-minimal graph on at least $N$ vertices
has a subgraph isomorphic to a subdivision of $M_{k}$.

\proof
We may assume that $k\ge65$.
Let $N$ be such that~\refthm{largenonplanar} holds for $k$ and $t:=130$,
and let $G$ be an X-minimal graph on at least $N$ vertices.
By~\refthm{shallow} the graph $G$ is $130$-shallow.
By~\refthm{largenonplanar} it has a subgraph isomorphic to a subdivision
of one of $B_k$, $M_k$, $K_{4,k}$, and $K'_{4,k}$.
But $G$ clearly has no subgraph isomorphic to a subdivision of 
$K_{4,k}$ or $K'_{4,k}$ (because the crossing number of these graphs is too 
large), and it has no subgraph isomorphic to a subdivision of
$B_k$ by~\refthm{noB}, because $k\ge65$.
Thus $G$ has a subgraph isomorphic to a subdivision of $M_k$, as desired.~\qed



\beginsection References

\def\JCTB{{\it J.\ Combin.\ Theory Ser.\ B}}



\bibitem{BokOpoRicSal} D.~Bokal, R.~B.~Richter and G.~Salazar,
Characterization of $2$-crossing-critical graphs I: Low connectivity or no
$V_{2n}$ minor,
manuscript, October 17, 2009. 

\bibitem{DinOpoOxlVer} G.~Ding, B.~Oporowski, J.~Oxley and D.~Vertigan,
Unavoidable minors of large $3$-connected matroids,
{\it\JCTB \bf 71} (1997), 244--299.

\bibitem{DinOpoThoVer}  G.~Ding, B.~Oporowski, R.~Thomas and D.~Vertigan,
Large $4$-connected nonplanar graphs,
manuscript, August 1999.

\bibitem{MohTho} B.~Mohar and C.~Thomassen, Graphs on surfaces,
Johns Hopkins University Press, Baltimore, MD, 2001.

\bibitem{OpoOxlTho}B.~Oporowski, J.~Oxley and R.~Thomas,
Typical subgraphs of 3- and 4-connected graphs,
\JCTB\ {\bf 57} (1993), 239--257.

\bibitem{RobSeyGM9}N.~Robertson and P.~D.~Seymour,
Graph Minors IX. Disjoint crossed paths,
\JCTB\ {\bf 49} (1990), 40--77.

\bibitem{RobSeyThoExt} N.~Robertson, P.~D.~Seymour and R.~Thomas,
Non-planar extensions of planar graphs,
manuscript.

\baselineskip 11pt
\vfill
\noindent
This material is based upon work supported by the National Science Foundation.
Any opinions, findings, and conclusions or
recommendations expressed in this material are those of the authors and do
not necessarily reflect the views of the National Science Foundation.

\end